\documentclass{article}
\usepackage[utf8]{inputenc}
\usepackage{amsthm}
\usepackage{amsmath}
\usepackage{authblk}

\title{Valid Widgets Contain Legal Subwidgets}
\author{Nathan Donagi}
\affil{Lower Merion High School}
\date{June 2022}

\begin{document}

\maketitle

\section{Introduction}

The purpose of this paper is to prove a linear algebra result. This has to do with the geometry of  ``widgets''. For us a widget is a collection of n pairs of points in a vector space. (The pairs represent the different possible spin states of a particle.) We investigate linear relations among such collections. A corollary of our theorem was conjectured in \cite{DO22} where it arose in an attempt to understand some issues in super string theory. In that paper an investigation of perturbative superstring theory with Ramond punctures required the special case when the ambient dimension is n. Here we prove the general case.

\smallskip
 
\section{Notation}

\smallskip

\begin{itemize}

    \item A widget is a set of n pairs of points in a vector space.

    \item A section is a subset containing at most one point from each pair of points.

    \item A legal widget is one in which every single section spans a space of dimension at most n-1.

    \item A subwidget is a subset of $k < n$ pairs out of the n pairs of points in a widget. (Note that this must be a proper subset.)

    \item A full widget is one in which the points span a space of dimension greater than or equal to n.

    \item A valid widget is one that is both legal and full.
    
    \item A point's opposite is the other point within a pair.

    \item Label the pairs $p_1$, $p_2$, $p_3$  etc.

    \item Label the points $p_1^+$, $p_1^-$, $p_2^+$, $p_2^-$ etc.
\end{itemize}

\section{Results}

\newtheorem{theorem}{Theorem}
\newtheorem{corollary}{Corollary}[theorem]
\newtheorem{lemma}[theorem]{Lemma}
    
Our main result is:
\begin{theorem}
Every valid widget contains a legal subwidget.
\end{theorem}
        
\begin{corollary}
A  widget is valid $\iff$ it is full and contains a legal subwidget.
\end{corollary}

\smallskip

It is fairly simple to see that if a widget contains a legal subwidget, the widget must be legal because choosing one point from each of the k pairs that span the subwidget will span at most a k-1 dimensional space. Therefore adding n-k extra points always span at most n-1 dimensions. 

\begin{corollary}

\label{nadias}
If a widget is valid then it contains a subwidget which is legal but not full. 

\end{corollary}

\smallskip

A valid widget must contain a legal subwidget, that subwidget is then either full or not. If it is full it is then also a valid widget and therefore must also contain a legal subwidget. Therefore by induction it is clear that one of the subwidgets of the original widget must therefore be legal but not full.

\smallskip

Corollary \ref{nadias} is a rephrasing of a conjecture in \cite{DO22}.

\section{Proof}

\begin{proof}

We denote the positive section: $g = {p_1^+, ..., p_{n}^+}$. 

\smallskip

For every i from 1 to n we let $g_i$ denote the section obtained by removing $p_i^+$ from g, and we let $s_i$ denotes the subwidget obtained by removing $p_i$ from the widget.

\smallskip

Assume there is a valid widget that does not contain a legal subwidget. Then no subwidget can be legal and therefore at least one section of the first n-1 pairs spans a n-1 space (otherwise the first n-1 pairs would be a legal subwidget). WLOG let that be $g_n$.

\smallskip

\begin{lemma}
Replacing $p_i^+$ by $p_i^+ + c * p_i^-$ for some constant c will keep a valid widget valid.
\end{lemma}

\smallskip

\begin{proof}

For all sections s taken from $s_i$, if s spans a n-1 dimensional space then for the widget to be legal both $p_i^+$ and $p_i^-$ must be in the span of s. Therefore adding some multiple of $p_i^-$ to $p_i^+$ will result in a new point that is also in the span of every section taken from $s_i$. Else since the span of s has at most dimension n-2, adding the point $p_i^+ + c * p_i^-$ will increase the dimension to at most n-1. Therefore the legality of the widget does not change. Also the fullness of the widget does not change as the total span of the widget has not changed.\end{proof}

\medskip

Call k the number of positive sections $g_i$ of a valid widget that span less than an n-1 space. We claim that for any valid widget that does not contain a legal subwidget, for some $k>0$, it is possible to create a new valid widget with a lower k. Let $g_i$ be one section that spans less than an n-1 space, with the span q of dimension n-2. The dimension of q can not be less than n-2 because at least one section $g_j$ spans a n-1 space.

\bigskip

If the opposites of every single point in $g_i$ are also in q, that means that every single point of each pair is in an n-2 space, and therefore there are n-1 pairs spanning a n-2 space, and therefore there is a legal subwidget.

\smallskip

So the opposite of some point in $g_i$ isn't in q. Therefore there is an infinite number of multiples c of that point’s opposite that can be added to it that then result in $g_i$ spanning a n-1 space. At most a finite number of these c could result in the span of a different section to go down in dimension, so there must therefore be a c that doesn't. Therefore it is possible to create a new widget with the same validity with a lower k.

\medskip

Through induction it is therefore possible to transform any widget that doesn't contain a legal subwidget into one in which either k=0 or it contains a legal subwidget. 

\medskip

If k=0 for the transformed widget, then every single section $g_i$ (for i from 1 to n) spans a n-1 space. Obviously then the sections g1 and g2 both span a n-1 space. The points $p_1^+$,$p_1^-$ must therefore be in the span of $g_1$. Therefore every point in $g_2$ is in the span of $g_1$ and since the span of both $g_2$ and $g_1$ have the same dimension they must be the same.  This exact same logic can be used to yield that the span of any section $g_i$ (for i from 1 to n) is the same n-1 space. Therefore every point in the transformed widget must be in that, so the transformed widget is therefore not full, and not valid. By Lemma 2 the original widget was not valid either, contradicting the assumption.

\end{proof}

\end{document}